\documentclass[11pt]{amsart}
\usepackage{amsopn}
\usepackage{amssymb, amscd}
\usepackage{amsmath}
\usepackage[table]{xcolor}

\usepackage{tikz}
\usetikzlibrary{matrix,positioning}
\usetikzlibrary{matrix,fit}
\usepackage{txfonts}
\usepackage{tikz}
\usetikzlibrary{matrix,backgrounds,fit}
\usepackage{mleftright}

\usetikzlibrary{matrix}

\usepackage{tikz}
\usetikzlibrary{decorations.pathreplacing}
\usetikzlibrary{decorations.text}

\newcommand{\nc}{\newcommand}

\nc{\eg}{\mathfrak{e} } \nc{\fg}{\mathfrak{f} } \nc{\vg}{\mathfrak{v} } \nc{\wg}{\mathfrak{w} }
\nc{\zg}{\mathfrak{z} } \nc{\ngo}{\mathfrak{n} } \nc{\kg}{\mathfrak{k} }
\nc{\mg}{\mathfrak{m} } \nc{\bg}{\mathfrak{b} } \nc{\ggo}{\mathfrak{g} }
\nc{\ggob}{\overline{\mathfrak{g}} } \nc{\sog}{\mathfrak{so} }
\nc{\sug}{\mathfrak{su} } \nc{\spg}{\mathfrak{sp} } \nc{\slg}{\mathfrak{sl} }
\nc{\glg}{\mathfrak{gl} } \nc{\cg}{\mathfrak{c} } \nc{\rg}{\mathfrak{r} }
\nc{\hg}{\mathfrak{h} } \nc{\tg}{\mathfrak{t} } \nc{\ug}{\mathfrak{u} }
\nc{\dg}{\mathfrak{d} } \nc{\ag}{\mathfrak{a} } \nc{\pg}{\mathfrak{p} }
\nc{\sg}{\mathfrak{s} } \nc{\affg}{\mathfrak{aff} }

\nc{\pca}{\mathcal{P}} \nc{\nca}{\mathcal{N}} \nc{\lca}{\mathcal{L}}
\nc{\oca}{\mathcal{O}} \nc{\mca}{\mathcal{M}} \nc{\tca}{\mathcal{T}}
\nc{\aca}{\mathcal{A}} \nc{\cca}{\mathcal{C}} \nc{\gca}{\mathcal{G}}
\nc{\sca}{\mathcal{S}} \nc{\hca}{\mathcal{H}} \nc{\bca}{\mathcal{B}}
\nc{\dca}{\mathcal{D}} \nc{\zca}{\mathcal{Z}}

\nc{\val}{\operatorname{val}}
\nc{\vp}{\varphi} \nc{\ddt}{\tfrac{d}{dt}} \nc{\dds}{\tfrac{{\rm d}}{{\rm d}s}}
\nc{\dpar}{\tfrac{\partial}{\partial t}} \nc{\im}{\mathtt{i}}

\nc{\SO}{\mathrm{SO}} \nc{\Spe}{\mathrm{Sp}} \nc{\Sl}{\mathrm{SL}}
\nc{\SU}{\mathrm{SU}} \nc{\Or}{\mathrm{O}} \nc{\U}{\mathrm{U}} \nc{\Gl}{\mathrm{GL}}
\nc{\Se}{\mathrm{S}} \nc{\Cl}{\mathrm{Cl}} \nc{\Spein}{\mathrm{Spin}}
\nc{\Pin}{\mathrm{Pin}} \nc{\G}{\mathrm{GL}_n(\RR)} \nc{\g}{\mathfrak{gl}_n(\RR)}

\nc{\RR}{{\Bbb R}} \nc{\HH}{{\Bbb H}} \nc{\CC}{{\Bbb C}} \nc{\ZZ}{{\Bbb Z}}
\nc{\FF}{{\Bbb F}} \nc{\NN}{{\Bbb N}} \nc{\QQ}{{\Bbb Q}} \nc{\PP}{{\Bbb P}}

\nc{\vs}{\vspace{.2cm}} \nc{\vsp}{\vspace{1cm}} \nc{\ip}{\langle\cdot,\cdot\rangle}
\nc{\ipp}{(\cdot,\cdot)} \nc{\la}{\langle} \nc{\ra}{\rangle} \nc{\unm}{\tfrac{1}{2}}
\nc{\unc}{\tfrac{1}{4}} \nc{\und}{\tfrac{1}{16}} \nc{\no}{\vs\noindent}
\nc{\lamkn}{\Lambda^2(\RR^{q+n})^*\otimes\RR^{q+n}} \nc{\lamn}{\Lambda^2(\RR^n)^*\otimes\RR^n} \nc{\lamp}{\Lambda^2\pg^*\otimes\pg}
\nc{\lamg}{\Lambda^2\ggo^*\otimes\ggo} \nc{\lamngo}{\Lambda^2\ngo^*\otimes\ngo}
\nc{\tangz}{{\rm T}^{\rm Zar}} \nc{\mum}{/\!\!/} \nc{\kir}{/\!\!/\!\!/}
\nc{\Ri}{\tfrac{4\Ricci_{\mu}}{||\mu||^2}} \nc{\ds}{\displaystyle}
\nc{\ben}{\begin{enumerate}} \nc{\een}{\end{enumerate}} \nc{\f}{\frac}
\nc{\lb}{[\cdot,\cdot]} \nc{\isn}{\tfrac{1}{||v||^2}}
\nc{\gkp}{(\ggo=\kg\oplus\pg,\ip)} \nc{\ukh}{(\ug=\kg\oplus\hg,\ip)}

\nc{\Hess}{\operatorname{Hess}} \nc{\ad}{\operatorname{ad}}
\nc{\Ad}{\operatorname{Ad}} \nc{\rank}{\operatorname{rank}}
\nc{\Irr}{\operatorname{Irr}} \nc{\End}{\operatorname{End}}
\nc{\Aut}{\operatorname{Aut}} \nc{\Inn}{\operatorname{Inn}}
\nc{\Der}{\operatorname{Der}} \nc{\Ker}{\operatorname{Ker}}
\nc{\Iso}{\operatorname{I}} \nc{\Diff}{\operatorname{Diff}}
\nc{\Lie}{\operatorname{Lie}} \nc{\tr}{\operatorname{tr}} \nc{\dif}{\operatorname{d}}
\nc{\sen}{\operatorname{sen}} \nc{\modu}{\operatorname{mod}}
\nc{\Riem}{\operatorname{Rm}} \nc{\Ricci}{\operatorname{Ric}}
\nc{\sym}{\operatorname{sym}} \nc{\symac}{\operatorname{sym^{ac}}}
\nc{\symc}{\operatorname{sym^{c}}} \nc{\scalar}{\operatorname{R}}
\nc{\grad}{\operatorname{grad}} \nc{\ricci}{\operatorname{Rc}}
\nc{\nr}{\operatorname{nr}} \nc{\riccic}{\operatorname{ric^{c}}}
\nc{\riccig}{\operatorname{ric^{\gamma}}} \nc{\Rin}{\operatorname{M}}
\nc{\Le}{\operatorname{L}} \nc{\tang}{\operatorname{T}}
\nc{\level}{\operatorname{level}} \nc{\rad}{\operatorname{r}}
\nc{\abel}{\operatorname{ab}} \nc{\CH}{\operatorname{CH}}
\nc{\mcc}{\operatorname{mcc}} \nc{\Adj}{\operatorname{Adj}}
\nc{\Order}{\operatorname{O}} \nc{\mm}{\operatorname{M}}
\nc{\inj}{\operatorname{inj}}  
\nc{\vol}{\operatorname{vol}} \nc{\Diag}{\operatorname{Diag}}

\theoremstyle{plain}
\newtheorem{theorem}{Theorem}[section]

\newtheorem{conjecture}[theorem]{Conjecture}

\theoremstyle{definition}
\newtheorem{definition}[theorem]{Definition}

\theoremstyle{remark}
\newtheorem{remark}[theorem]{Remark}

\newtheorem{example}[theorem]{Example}

 \nc{\Pf}{\operatorname{\textsf{Pf}\,}}

\title[Pfaffian representations of plane curves]{Pfaffian representations of plane curves}
\author{David Oscari}
\address{FaMAF and CIEM, Universidad Nacional de C\'ordoba, C\'ordoba, Argentina}
\email{oscari@famaf.unc.edu.ar}
\thanks{This research was partially supported by grants from CONICET, FonCyT (Argentina) and SeCyT (Universidad Nacional de C\'ordoba)}

\begin{document}


\maketitle

\begin{abstract}
Let $R$ be a commutative ring with 1. For every homogeneous polynomial $f(X_0,X_1,X_2)$ in $R[X_0,X_1,X_2]$ of degree $d\leq25$, we find a {\it explicit} linear Pfaffian $R$-representation of $f$. We describe an empirical method that leads us to find such $R$-representations. This generalizes  and constitutes an alternative proof (up to degree 25) of a result due to A. Beauville \cite{Bea} about the existence of linear Pfaffian $\mathbb{K}$-representations for any smooth plane curve of degree $d\geq2$, where $\mathbb{K}$ is an algebraically closed field of characteristic zero.
\end{abstract}

\section{Introduction}

Let $R$ be a commutative ring with $1$, and let $M=[m_{ij}]$ be skew-symmetric matrix. There are several equivalent definitions of the Pfaffian  of a skew-symmetric matrix. In our computations with Maple\textsuperscript{\texttrademark}, we have used a Laplace-type expansion along last row. For other definitions of Pfaffian, its pro\-per\-ties and applications we refer the reader to \cite{H}, \cite{DW}, \cite[Chapter 7]{G}, \cite[Appendix D]{FP}. For a history of the Pfaffians we recommend \cite{K}.
\begin{definition}[Expansion along last row]\label{definicion recursiva}
The  \textit{Pfaffian} $\Pf(M)$ of a $2d\times2d$ skew-symmetric matrix $M=[m_{ij}]$ with entries in $R$ is defined recursively as
\begin{align*}
 \Pf\left(\left(
  \begin{array}{cc}
    0 & m_{12} \\
    -m_{12} & 0
    \end{array}
  \right)\right) & :=   m_{12}  & \textrm{ if }d=1\,, \\
\Pf(M)   & :=  \sum_{j=1}^{2d-1}(-1)^j\cdot m_{1j}\cdot \Pf\!\left(M^{\{2d,j\}}\right)   &  \textrm{ if }d\geq2\,,
\end{align*}
where $M^{\{2d,j\}}$ is the $(2d-2)\times(2d-2)$-matrix obtained from $M$ by deleting both the $(2d,j)$-th rows and $(2d,j)$-th columns.
\end{definition}


\begin{definition}\label{representation pfaffiana linear} Let $R$ be a commutative ring with $1$ and let $f(X_0,\dots,X_n)$ be a homogeneous polynomial in
$R[X_0,\dots,X_n]$ of degree $d$. Then $f$ admits a {\it linear Pfaffian $R$-representation} if there exist $2d\times2d$ skew-symmetric matrices $A_0,\dots,A_n$ with
entries in $R$ such that
$$
\Pf(X_0A_0+\cdots+X_nA_n)= f(X_0,\dots,X_n).
$$
\end{definition}
\noindent Such a matrix $M=X_0A_0+\cdots+X_nA_n$ is said to be a linear Pfaffian $R$-representation of $f$.

In the case $n=2$, Beauville showed, among others, that every smooth plane curve $f(X_0,X_1,X_2)$ in $\mathbb{K}[X_0,X_1,X_2]$ of degree $d\geq 2$ admits a linear Pfaffian $\mathbb{K}$-representation, where  $\mathbb{K}$ is an algebraically closed field of characteristic zero, see \cite[Proposition 5.1]{Bea}. Buckley and K\v{o}sir \cite{BK} parametrise all linear Pfaffian representations of a plane curve. Two Pfaffian $\mathbb{K}$-representations
$M$ and $M'$ are equivalent if there exists $X\in\text{GL}_{2d}(\mathbb{K})$ such that $M'=XMX^t$.
In \cite{B} the elementary transformations of linear Pfaffian representations are considered, proving that every two Pfaffian representations of a plane curve
$f(X_0,X_1,X_2)$ of degree $d$ can be bridged by a finite sequence of elementary transformations.

In the case $n=3$, Beauville showed that every smooth cubic surface admits a linear Pfaffian representation \cite[Proposition 7.6(a)]{Bea}. Then Fania and Mezzetti \cite{FM} showed that every cubic surface admits such a linear Pfaffian representation (see also \cite{FMb}).

The fact that a {\it general} surface of degree $d$ admits a linear Pfaffian  representation if and only if $d\leq15$ had been proved by Beauville-Schreyer \cite[Proposition 7.6(b)]{Bea}. For $d=4$, Coskun, Kulkarni and Mustopa \cite{CKM} generalize  the previous result: every smooth quartic surface admits a linear Pfaffian representation. Faenzi \cite{F} proved that a general surface of degree $d$ admits a almost quadratic Pfaffian representation if and only if $d\leq15$.
In \cite{CF} the Pfaffian representations with entries almost linear of general surfaces  are considered.

In \cite[Theorem (47.3)]{AR} Adler proved that a {\it generic} cubic threefold ($n=4$) admits a Pfaffian representation. Then Beauville proved that every smooth cubic threefold admits a Pfaffian representation \cite[Proposition 8.5]{Bea}. On the other hand, Iliev and Markushevich showed that a {\it generic} quartic threefold admits a Pfa\-ffian representation \cite[Proposition 1.2]{IM}. Brambilla and Faenzi \cite{BF} also provide a proof, which is not computer-aided, about existence of Pfaffian representations for quartic threefold. The method used in both \cite{AR} and \cite{IM} is to prove that the differential of the Pfa\-ffian map is of maximal rank in a particular skew-symmetric matrix $M_0$ with linear entries.


Independently, F. Tanturri and F. Han were interested on the explicit construction of linear Pfaffian representations for cubic surfaces.
In \cite{T} and \cite{Han} algorithms are provided whose inputs involve a cubic surface $f(x,y,z,t)$ in $\mathbb{K}[x,y,z,t]$  and whose output is a linear Pfaffian representation of $f$.

In this note, we focus on homogeneous polynomials $f(x,y,z)$ in $R[x,y,z]$ of degree $d\geq2$.

Our starting point was one Dickson's idea modified \cite{D}. It was difficult to get explicit linear  pfa\-ffian $R$-representations of $f(x,y,z)$ for degrees $d\leq7$. In acordding to our observations, we propose the following.

\begin{conjecture}\label{main}
Let $R$ be a commutative ring with 1.
\begin{enumerate}
    \item Every homogeneous polynomial $f(X_0,X_1,X_2)$ in $R[X_0,X_1,X_2]$ of degree $d\geq2$ admits a linear Paffian $R$-representation.
    \item Such a $R$-representation can be find it by using the empirical method on Section \ref{empirical method}, $d\geq 5$.
\end{enumerate}
\end{conjecture}

\begin{remark}
By applying that empirical method, we have checked the Conjecture \ref{main} for all degree $5\leq d\leq25$. For $d\leq5$, in \cite{O} we have obtained {\it explicit} $R$-representatios by a less efficient method.
\end{remark}


\section{Empirical Method} \label{empirical method}

\newcommand{\X}{x}
\newcommand{\Y}{y}
\newcommand{\Z}{z}
Fix $f(\X,\Y,\Z)$ in $R[\X,\Y,\Z]$ of degree $d\geq5$:
\begin{align*}
   f(\X,\Y,\Z) & = \sum_{\substack{0\,\leq\,i,j,k\, \leq\, d \vspace{0.1em}\\ i+j+k=d}}\: \Theta_{ijk} \cdot x^iy^jz^k  \\
               &=x^d\cdot\Theta_{d00}+y^d\cdot\Theta_{0d0}+z^d\cdot\Theta_{00d}+  \sum_{\substack{0\,\leq\,i,j,k\, \leq\, d-1 \vspace{0.1em}\\ i+j+k=d}}\: \Theta_{ijk} \cdot x^iy^jz^k.\\
\end{align*}

Let we consider the $2d\times2d$ skew-symmetric matrix
$$
\newcommand{\xpaso}{0}  
\newcommand{\ypaso}{0}  
\newcommand{\xcell}{1}
\newcommand{\ycell}{1}
M=
\left[\begin{array}{c}
    {
    \footnotesize
 \begin{tikzpicture}[scale=.8]
 \node  at (0 + \xpaso + \xcell *   3, 0 + \ypaso - \ycell  *  7) {$*$};
 \node  at (0 + \xpaso + \xcell *   5, 0 + \ypaso - \ycell  *  12) {$*$};
 \node  at (0 + \xpaso + \xcell *   1, 0 + \ypaso - \ycell  *  1) {$0$};
 \node  at (0 + \xpaso + \xcell *   2, 0 + \ypaso - \ycell  *  1) {$\textcolor{black}{x\Theta_{d00}}$};
 \node  at (0 + \xpaso + \xcell *   3, 0 + \ypaso - \ycell  *  1) {$m_{1,3}$};
 \node  at (0 + \xpaso + \xcell *   4, 0 + \ypaso - \ycell  *  1) {$m_{1,4}$};
 \node  at (0 + \xpaso + \xcell *   5, 0 + \ypaso - \ycell  *  1) {$m_{1,5}$};
 \node  at (0 + \xpaso + \xcell *   8, 0 + \ypaso - \ycell  *  1) {$m_{1,d+1}$};
 \node  at (0 + \xpaso + \xcell *   2, 0 + \ypaso - \ycell  *  2) {$0$};
 \node  at (0 + \xpaso + \xcell *   3, 0 + \ypaso - \ycell  *  2) {$\textcolor{black}{y\Theta_{0d0}}$};
 \node  at (0 + \xpaso + \xcell *   4, 0 + \ypaso - \ycell  *  2) {$m_{2,4}$};
 \node  at (0 + \xpaso + \xcell *   5, 0 + \ypaso - \ycell  *  2) {$m_{2,5}$};
 \node  at (0 + \xpaso + \xcell *   8, 0 + \ypaso - \ycell  *  2) {$m_{2,d+1}$};
 \node  at (0 + \xpaso + \xcell *   3, 0 + \ypaso - \ycell  *  3) {$0$};
 \node  at (0 + \xpaso + \xcell *   4, 0 + \ypaso - \ycell  *  3) {$\textcolor{black}{z\Theta_{00d}}$};
 \node  at (0 + \xpaso + \xcell *   5, 0 + \ypaso - \ycell  *  3) {$m_{3,5}$};
 \node  at (0 + \xpaso + \xcell *   8, 0 + \ypaso - \ycell  *  3) {$m_{3,d+1}$};
 \node  at (0 + \xpaso + \xcell *   4, 0 + \ypaso - \ycell  *  4) {$0$};
 \node  at (0 + \xpaso + \xcell *   5, 0 + \ypaso - \ycell  *  4) {$m_{4,5}$};
 \node  at (0 + \xpaso + \xcell *   8, 0 + \ypaso - \ycell  *  4) {$m_{4,d+1}$};
 \node  at (0 + \xpaso + \xcell *   5, 0 + \ypaso - \ycell  *  5) {$0$};
 \node  at (0 + \xpaso + \xcell *   8, 0 + \ypaso - \ycell  *  5) {$m_{5,d+1}$};
 \node  at (0 + \xpaso + \xcell *   8, 0 + \ypaso - \ycell  *  7) {$m_{d,d+1}$};
 \node  at (0 + \xpaso + \xcell *   8, 0 + \ypaso - \ycell  *  8) {$0$};
\newcommand{\lineapunto}{.\ .\ .\ .}
\path[decorate,decoration={text along path,text={\lineapunto},text align=center}]
(0 + \xpaso + \xcell *   6, 0 + \ypaso - \ycell  *  1) -- (0 + \xpaso + \xcell *   7, 0 + \ypaso - \ycell  *  1);
\path[decorate,decoration={text along path,text={\lineapunto},text align=center}]
(0 + \xpaso + \xcell *   6, 0 + \ypaso - \ycell  *  2) -- (0 + \xpaso + \xcell *   7, 0 + \ypaso - \ycell  *  2);
\path[decorate,decoration={text along path,text={\lineapunto},text align=center}]
(0 + \xpaso + \xcell *   6, 0 + \ypaso - \ycell  *  3) -- (0 + \xpaso + \xcell *   7, 0 + \ypaso - \ycell  *  3);
\path[decorate,decoration={text along path,text={\lineapunto},text align=center}]
(0 + \xpaso + \xcell *   6, 0 + \ypaso - \ycell  *  4) -- (0 + \xpaso + \xcell *   7, 0 + \ypaso - \ycell  *  4);
\path[decorate,decoration={text along path,text={\lineapunto},text align=center}]
(0 + \xpaso + \xcell *   6, 0 + \ypaso - \ycell  *  5) -- (0 + \xpaso + \xcell *   7, 0 + \ypaso - \ycell  *  5);
\path[decorate,decoration={text along path,text={\lineapunto},text align=center}]
(0 + \xpaso + \xcell *   6, 0 + \ypaso - \ycell  *  6) -- (0 + \xpaso + \xcell *   7, 0 + \ypaso - \ycell  *  7);
  \path[decorate,decoration={text along path,text={\lineapunto},text align=center}]
(0 + \xpaso + \xcell *   8, 0 + \ypaso - \ycell  *  6.5) -- (0 + \xpaso + \xcell *   8, 0 + \ypaso - \ycell  *  5.5);
\newcommand{\extra}{0}
\newcommand{\ladoalado}{0}
\newcommand{\arribaabajo}{0}
\draw   (9+\ladoalado,-.5+\extra) -- (9+\ladoalado,-14.5-\extra);   
\draw   (11+\ladoalado,-.5+\extra) -- (11+\ladoalado,-2.5-\extra);   
\draw   (9+\ladoalado,-2.5+\arribaabajo) -- (15+\ladoalado,-2.5+\arribaabajo-\extra);     
\draw   (1+\ladoalado,-8.5+\arribaabajo) -- (15+\ladoalado,-8.5+\arribaabajo-\extra);     
%
\newcommand{\xnil}{.5}
\newcommand{\ynil}{0}
 \node  at (0 + \xnil + \xcell *   9,  0 + \ynil - \ycell  *  1) {$-y$};
 \node  at (0 + \xnil + \xcell *  10,  0 + \ynil - \ycell  *  1) {$-z$};
 \node  at (0 + \xnil + \xcell *   9,  0 + \ynil - \ycell  *  2) {$-z$};
 \node  at (0 + \xnil + \xcell *   10,  0 + \ynil - \ycell  *  2) {$0$};
 \node  at (0 + \xnil + \xcell *   9,  0 + \ynil - \ycell  *  8) {$x$};
 \node  at (0 + \xnil + \xcell *   10,  0 + \ynil - \ycell  * 7) {$x$};
 \node  at (0 + \xnil + \xcell *   10,  0 + \ynil - \ycell  * 8) {$y$};
 \node  at (0 + \xnil + \xcell *   11,  0 + \ynil - \ycell  * 6) {$x$};
 \node  at (0 + \xnil + \xcell *   11,  0 + \ynil - \ycell  * 7) {$-y$};
 \node  at (0 + \xnil + \xcell *   11,  0 + \ynil - \ycell  * 8) {$(-1)^dz$};
 \node  at (0 + \xnil + \xcell *   12,  0 + \ynil - \ycell  * 5) {$x$};
 \node  at (0 + \xnil + \xcell *   12,  0 + \ynil - \ycell  * 6) {$-y$};
 \node  at (0 + \xnil + \xcell *   12,  0 + \ynil - \ycell  * 7) {$-z$};
  \node  at (0 + \xnil + \xcell *   14,  0 + \ynil - \ycell  * 3) {$x$};
 \node  at (0 + \xnil + \xcell *   14,  0 + \ynil - \ycell  * 4) {$-y$};
 \node  at (0 + \xnil + \xcell *   14,  0 + \ynil - \ycell  * 5) {$-z$};
\path[decorate,decoration={text along path,text={\lineapunto},text align=center}]
(0 + \xnil + \xcell *   13.5, 0 + \ynil - \ycell  * 3.5) -- (0 + \xnil + \xcell *   12.5, 0 + \ynil - \ycell  *  4.5);
\path[decorate,decoration={text along path,text={\lineapunto},text align=center}]
(0 + \xnil + \xcell *   13.5, 0 + \ynil - \ycell  * 4.5) -- (0 + \xnil + \xcell *   12.5, 0 + \ynil - \ycell  *  5.5);
\path[decorate,decoration={text along path,text={\lineapunto},text align=center}]
(0 + \xnil + \xcell *   13.5, 0 + \ynil - \ycell  * 5.5) -- (0 + \xnil + \xcell *   12.5, 0 + \ynil - \ycell  *  6.5);
\node  at (0 + \xnil + \xcell *   10.5,  0 + \ynil - \ycell  * 4) { \large$\mathbf{0}$};
\node  at (0 + \xnil + \xcell *   12.5,  0 + \ynil - \ycell  * 1.5) { \large$\mathbf{0}$};
\node  at (0 + \xnil + \xcell *   14,  0 + \ynil - \ycell  * 7.5) { \large$\mathbf{0}$};
\node  at (0 + \xnil + \xcell *   12,  0 + \ynil - \ycell  * 11.5) { \large$\mathbf{0}_{d-1}$};
%
%
\newcommand{\xm}{-.5}
\newcommand{\ym}{0}
\end{tikzpicture}
}
\end{array}\right]
$$

\noindent where $m_{ij}=a_{ij}x+b_{ij}y+c_{ij}z$. For all degree $5\leq d\leq 25$, we have checked that holds:
\begin{itemize}
    \item
$\displaystyle
\Pf(M)=x^d\cdot\Theta_{d00}+y^d\cdot\Theta_{0d0}+z^d\cdot\Theta_{00d}+
\sum_{\substack{0\,\leq\,i,j,k\, \leq\, d-1\vspace{0.1em}\\ i+j+k=d }}\: e_{ijk} \cdot x^iy^jz^k.
$
\item The $e_{ijk}$' s are polynomials of degree 1 in the indeterminates $\Theta_{ijk},a_{ij}, b_{ij}, c_{ij}$.
\item Since $\Pf(M)=f(x,y,z)$, we obtain a linear system of equations:
      $$
      \text{(S)}\qquad \Big\{  e_{ijk}=\Theta_{ijk} \,, \qquad\qquad  i+j+k=d \text{ and }0\leq\,i,j,k\, \leq\, d-1
      $$
where the $\Theta_{ijk}$'s are parameters.
\item The linear system of equations (S) always has a solution in $R$.
\end{itemize}
\noindent  Therefore, $f(x,y,z)$ admits a explicit linear Pfaffian $R$-representation for all degree $5\leq d\leq25$.

\begin{example}\label{ejemplo}
Consider an arbitrary homogeneous polynomial $f(x,y,z)$ in $R[x,y,z]$ of degree $d=5$. Instead of the notation $\Theta_{ijk}$, we simply enumerate their coefficients, e.g.
\begin{align*}
f(x,y,z)&= x^{5}\Theta_1+y^{5}\Theta_2+x^{5}\Theta_3\,  + \\
        & \qquad +x^{4}y\Theta_{4}+x^{3}y^{2}\Theta_{5}+x^{2}y^{3}\Theta_{6}+xy^{4}\Theta_{7}+x^{4}z\Theta_{8}+x^{3}z^{2}\Theta_{9}+x^{2}z^{3}\Theta_{10}+xz^{4}\Theta_{11}  \\
        & \qquad  +y^{4}z\Theta_{12} +y^{3}z^{2}\Theta_{13}+y^{2}z^{3}\Theta_{14}+yz^{4}\Theta_{15}+xyz^{3}\Theta_{16}+xy^{2}z^{2}\Theta_{17}+xy^{3}z\Theta_{18}   \\
        & \qquad  +x^{2}yz^{2}\Theta_{19}+x^{2}y^{2}z\Theta_{20}  +x^{3}yz\Theta_{21}  \,.     \\
\end{align*}
Let we consider the $10\times10$ skew-symmetric matrix
{\small
$$
M =  \left[\: \begin {array}{cccccccccc}
0& \textcolor{black}{x\Theta_1} & m_{2} &m_{3}&m_{4}&m_{5}&\multicolumn{1}{|c}{-y}&\multicolumn{1}{c|}{-z} & 0& 0\\
   &0& \textcolor{black}{y\Theta_2} & m_{11} &m_{12}&m_{13}&\multicolumn{1}{|c}{-z}&\multicolumn{1}{c|}{0}&0 &0\\
  \cline{7-10}
    &&0& \textcolor{black}{y\Theta_3} & m_{19} &m_{20}&\multicolumn{1}{|c}{0} &0 & 0&x\\
  & & &0&m_{25}&m_{26}&\multicolumn{1}{|c}{ 0}& 0&x&  -y \\
  & * &  & &0&m_{31}&\multicolumn{1}{|c}{0 }&x&-y& -z \\
  & & & & &0&\multicolumn{1}{|c}{x}&y&-z&0\\
\cline{1-10}
 & & & & & &\multicolumn{1}{|c}{0 }&0&0&0\\
 & & & & & &\multicolumn{1}{|c}{ 0}&0&0&0\\
 & & & * & & &\multicolumn{1}{|c}{ 0}&0&0&0\\
 & & & & & &\multicolumn{1}{|c}{ 0}&0&0&0\\
\end {array} \: \right]
$$
}
where $m_{r}=a_{r}x+b_{r}y+c_{r}z$. In this case
\begin{align*}
& \Pf(M)=  \\
   & x^5\Theta_1+y^5\Theta_2+z^5\Theta_3+ {x}^{4}y\cdot \left( b_{{1}}-a_{{13}} \right) +  {x}^{3}{y}^{2}\cdot \left( a_{{12}}-b_{{13}} \right)+  {x}^{2}{y}^{3}\cdot \left( a_{{11}}+b_{{12}} \right) +\\
   &\qquad +  x{y}^{4}\cdot \left( a_{{10}}+b_{{11}} \right) + {x}^{4}z\cdot \left( -a_{{12}}+c_{{1}}+a_{{5}} \right) +
    {x}^{3}{z}^{2}\cdot \left( a_{{3}}-a_{{10}}-c_{{12}}+a_{{31}}+c_{{5}} \right) +\\
   &\qquad +{x}^{2}{z}^{3}\cdot \left( c_{{3}}-c_{{10}}+a_{{20}}-a_{{25}}+c_{{31}} \right) +
   x{z}^{4}\cdot \left( c_{{20}}+a_{{18}}-c_{{25}} \right) +
  z{y}^{4}\cdot \left( -b_{{2}}+c_{{10}} \right)+  \\
   &\qquad + {z}^{2}{y}^{3}\cdot \left( b_{{20}}-c_{{2}} \right)+
  {z}^{3}{y}^{2}\cdot \left( c_{{20}}-b_{{19}} \right) +  {z}^{4}y\cdot \left( -c_{{19}}+b_{{18}} \right)+ \\
   & \qquad +
  xy{z}^{3}\cdot \left( -a_{{19}}+c_{{26}}+b_{{20}}-b_{{25}} \right) +
 x{y}^{2}{z}^{2}\cdot  \left( -c_{{10}}-c_{{3}}+b_{{26}}+a_{{20}} \right) + \\
   &\qquad +
  x{y}^{3}z\cdot \left( -b_{{3}}+c_{{11}}-\Theta_{{2}}-a_{{2}} \right)+
  {x}^{2}y{z}^{2}\cdot \left( -\Theta_{{2}}-2\,c_{{11}}+b_{{3}}+a_{{26}}+b_{{31}}-c_{{4}} \right)+\\
   &\qquad +
  {x}^{2}{y}^{2}z\cdot \left( -2\,b_{{11}}+c_{{12}}-a_{{10}}-a_{{3}}-b_{{4}} \right)+\\
   & \qquad +
 {x}^{3}yz\cdot  \left( -c_{{13}}+b_{{5}}-a_{{4}}-b_{{12}}-2\,a_{{11}} \right)\,.  \\
\end{align*}
Since $\Pf(M)=f(x,y,z)$, by equaliting the correnponding coefficients, we obtain a linear system of equations where the $\Theta_i$'s are parameters:
{\small
\begin{equation}\tag{$\text{S}'$}\label{sprima}
\left\{
 \begin{array}{rccrrrc}
b_{{1}}-a_{{13}}&=\Theta_{{4}}  &   \qquad \qquad     &     -c_{{2}}+b_{{20}}&=\Theta_{{13}}     &   \\
-b_{{13}}+a_{{12}}&=\Theta_{{5}} &         &   c_{{20}}-b_{{19}}&=\Theta_{{14}}  &   \\
a_{{11}}+b_{{12}}&=\Theta_{{6}}  &  &  -c_{{19}}+b_{{18}}&=\Theta_{{15}}      &   \\
b_{{11}}+a_{{10}}&=\Theta_{{7}} &   &  c_{{26}}-b_{{25}}-a_{{19}}+b_{{20}}&=\Theta_{{16}}         &   \\
c_{{1}}-a_{{12}}+a_{{5}}&=\Theta_{{8}} &   &  -c_{{3}}+b_{{26}}+a_{{20}}-c_{{10}}&=\Theta_{{17}}         &   \\
a_{{31}}-c_{{12}}+c_{{5}}+a_{{3}}-a_{{10}}&=\Theta_{{9}} &     &  -\Theta_{{2}}+c_{{11}}-b_{{3}}-a_{{2}}&=\Theta_{{18}}        &   \\
c_{{31}}-a_{{25}}+c_{{3}}+a_{{20}}-c_{{10}}&=\Theta_{{10}}  &  & -\Theta_{{2}}+b_{{31}}-2\,c_{{11}}-c_{{4}}+a_{{26}}+b_{{3}}&=\Theta_{{19}}  &        \\
-c_{{25}}+c_{{20}}+a_{{18}}&=\Theta_{{11}}  &  &  -2\,b_{{11}}-b_{{4}}-a_{{3}}+c_{{12}}-a_{{10}}&=\Theta_{{20}}  &         \\
-b_{{2}}+c_{{10}}&=\Theta_{{12}} &   & -2\,a_{{11}}-a_{{4}}-c_{{13}}-b_{{12}}+b_{{5}}&=\Theta_{{21}}  &          \\
 \end{array}
 \right.
\end{equation}
}

\noindent We can solve linear system (\ref{sprima}) by using e.g. Maple\textsuperscript{\texttrademark}. The general solution has 42 components, where there are 24 free parameters, say $t_1,\dots, t_{24}$. Each of their components is a polynomial in $\Theta_1,\dots,\Theta_{21}$ and $t_1,\dots, t_{24}$ (with coefficients in $R$), i.e. each component of the general solution is an element in $R$.

\end{example}

\bibliographystyle{plain}

\begin{thebibliography}{10}


\bibitem[AR]{AR}
{\sc A.~Adler and S.~Ramanan.}
\newblock {\em Moduli of abelian varieties}, Lecture Notes in Mathematics {\bf 1644}.
\newblock Springer-Verlag, Berlin, 1996.


\bibitem[Bea]{Bea}
\textsc{A. Beauville}.
\newblock Determinantal hypersurfaces.
\newblock {\it Mich. Math. J.} {\bf 48}(1): 39--64, 2000.


\bibitem[BF]{BF}
{\sc M.~C. Brambilla and D.~Faenzi.}
\newblock Moduli spaces of rank-2 {ACM} bundles on prime {F}ano threefolds.
\newblock {\em Mich. Math. J.} {\bf 60}(1): 113--148, 2011.


\bibitem[B]{B}
{\sc A. Buckley.}
\newblock Elementary transformations of Pfaffian representations of plane curves.
\newblock {\it Linear Algebra Appl.}  {\bf 443} (4): 758--780, 2010.


\bibitem[BK]{BK}
{\sc A. Buckley, T. Ko\v{s}ir.}
\newblock Plane Curves as Pfaffians.
\newblock {\it Sc. Norm. Super. Pisa Cl. Sci.} (5) {\bf 10} (2): 363--388, 2011.


\bibitem[CF]{CF}
{\sc L.~Chiantini, D.~Faenzi.}
\newblock On general surfaces defined by an almost linear {P}faffian.
\newblock {\it Geom. Dedicata} 142:91--107, 2009.


\bibitem[CKM]{CKM}
{\sc E.~Coskun, R.~S. Kulkarni, Y.~Mustopa.}
\newblock Pfaffian quartic surfaces and representations of Clifford algebras.
\newblock {\it Documenta Math.} {\bf 17}: 91--107 (2012).


\bibitem[D]{D}
{\sc L. E. Dickson}.
\newblock Determination of all general homogeneous polynomials expressible as determinants with linear elements.
\newblock {\it Trans. Amer. Math. Soc.} {\bf 22}: 167--179, 1921.


\bibitem[DW]{DW}
{\sc A. W. M. Dress, W. Wenzel.}
\newblock A simple proof of an identity concerning pfaffians of skew symmetric matrices.
\newblock {\it Adv. in Math.} {\bf 112}: 120--134, 1995.


\bibitem[F]{F}
{\sc D.~Faenzi.}
\newblock A remark on {P}faffian surfaces and a{CM} bundles.
\newblock In {\em Vector bundles and low codimensional subvarieties: state of
  the art and recent developments}, volume~21 of {\em Quad. Mat.}, pages
  209--217. Dept. Math., Seconda Univ. Napoli, Caserta, 2007.


\bibitem[FM]{FM}
{\sc M.~L. Fania, E.~Mezzetti.}
\newblock On the {H}ilbert scheme of {P}alatini threefolds.
\newblock {\em Adv. Geom.}, {\bf 2}(4): 371--389, 2002.


\bibitem[FMb]{FMb}
{\sc M.~L. Fania, E.~Mezzetti.}
\newblock Erratum to ``On the Hilbert scheme of Palatini threefolds''.
\newblock {\em Adv. Geom.} {\bf 8}(1): 153--154, 2008.


\bibitem[FP]{FP}
{\sc W. Fulton, P. Pragacz.}
\newblock {\it Shubert varieties and degeneracy loci}.
\newblock Lecture Notes in Mathematics {\bf  1689}, Springer-Verlag, 1998.


\bibitem[G]{G}
{\sc C.D.~Godsil.}
\newblock {\it Algebraic Combinatorics}.
\newblock Chapman and Hall, 1993.


\bibitem[H]{H}
{\sc A.M.~Hamel.}
\newblock Pfaffian Identities: A Combinatorial Approach.
\newblock {\it J. Combin. Theory Serie A} {\bf 94}: 205--217, 2001.


\bibitem[Han]{Han}
{\sc F.~Han.}
\newblock Pfaffian bundles on cubic surfaces and configurations of planes.
\newblock {\it Math. Z.} {\bf 278}(1): 363--383, 2014.


\bibitem[IM]{IM}
{\sc A.~Iliev, D.~Markushevich.}
\newblock Quartic 3-fold: {P}faffians, vector bundles, and half-canonical curves.
\newblock {\em Michigan Math. J.} {\bf 47}(2): 385--394, 2000.


\bibitem[K]{K}
{\sc D. E. Knuth}.
\newblock Overlapping Pfaffians.
\newblock {\it  Electron. J. Combin.} {\bf  3}: 1--13, 1996.


\bibitem[O]{O}
{\sc D.~Oscari.}
\newblock Explicit linear pfaffian representations of plane curves up to degree 5.
\newblock Preprint arXiv:1712.03600v1.


\bibitem[T]{T}
{\sc F. Tanturri}.
\newblock Paffian representation of cubic surfaces.
\newblock {\it Geom. Dedicata} {\bf 168}: 69--86, 2014.


\end{thebibliography}

\def\cprime{$'$}

\end{document}